\documentstyle[12pt,a4,amstex,twoside]{article}

\def\mytitle{On The Structure of\\ \ Quantum Super Groups  GL$_q(m|n)$}
\def\myemail{{\small\tt phung@@ictp.trieste.it}}
\def\myaddress{\small\it International Center for Theoretical Physics\\\small\it P.O. Box 586, 34100 Trieste, Italy}
\def\myperaddress{Hanoi Institute of Mathematics P.O. Box 631, 10000 Bo Ho, Hanoi, Vietnam}
\def\mythanks{The author would like to thank the International Center for Theoretical Physics for its hospitality and financial support.}
\def\myabstract{We show that a quantum super matrix in standard format is invertible if and only if its block matrices of even entries are invertible. We prove the q-analog of the well-known formula for the Berezinian.}

\def\part{\vdash}

\def\myauthor{Ph\`ung H{{\accent"5E o}\kern-.28em\raise.2ex\hbox{\char'22}\kern-.20em}
 H{a\kern-.370em\raise.16ex\hbox{\char'47}\kern.1em}i}

\def\iii{^{-1}}

\def\lora{\longrightarrow}
\def\ot{\otimes}

\def\loma{\longmapsto}

\def\si{\sigma}
\newcommand{\bbas}{\begin{eqnarray*}}
\newcommand{\eeas}{\end{eqnarray*}}
\newcommand{\bbar}{\begin{array}}
\newcommand{\eear}{\end{array}}
\newcommand{\bbs}{\begin{displaymath}}
\newcommand{\ees}{\end{displaymath}}
\newcommand{\bb}{\begin{equation}}
\def\ee{\end{equation}}
\def\eea{\end{eqnarray}}

\newcommand{\dayso}[1]{\mbox{$1,2,\ldots ,#1$}}

\newcommand{\daysub}[2]{\mbox{$#1_1,#1_2,\ldots ,#1_{#2}$}}
\def\bba{\begin{eqnarray}}
\newtheorem{thm}{Theorem}[section]
\newtheorem{lem}[thm]{Lemma}

\newtheorem{edl}[thm]{Theorem}

\newtheorem{cor}[thm]{Corollary}

\newtheorem{pro}[thm]{Proposition}

\def\Lambd{{\mit\Lambda}}

\def\det{\mbox{\rm det{}}}

\def\sdet{\mbox{\rm Ber$_q$}{}}
\newcommand{\Ai}{\mbox{$A^{-1}$}{}}
\newcommand{\Vi}{\mbox{$V^{-1}$}{}}
\newcommand{\Hi}{\mbox{$H^{-1}$}{}}
\newcommand{\Di}{\mbox{$D^{-1}$}{}}
\def\hi{\hat{i}}
\def\hj{\hat{j}}
\def\hm{\hat{m}}

\def\hn{\hat{n}}
\def\hk{\hat{k}}
\def\hl{\hat{l}}
\def\hp{\hat{p}}
\def\hq{\hat{q}}

\def\hs{\hat{s}}

\def\Im{\mbox{\rm Im}{}}

\def\Ber{\mbox{\rm Ber{}}}

\def\kk{{\Bbb K}}
\def\d{\mbox{\rm D{}}}
\def\H{{\cal H}}

\def\eqeee{\mbox{\rule{.75ex}{1.5ex}}}
\def\eee{\rule{.75ex}{1.5ex}\vskip1ex}

\def\proof{{\it Proof.\ }}
\newcommand{\va}{\varepsilon}
\def\nml{\normalsize}

\newcommand{\Tr}{\mbox{\rm Tr}}
\renewcommand{\dim}{\mbox{\rm dim}}

\def\rref#1{(\ref{#1})}

\font\Fraktur=eufm10 scaled\magstep1          
   \newcommand{\fraktur}[1]{\mbox{\Fraktur #1}}  %
   \font\Fraktu=eufm7 scaled\magstep1            
   \newcommand{\fraktu}[1]{\mbox{\Fraktu #1}}    %
   \font\Frakt=eufm5 scaled\magstep1             
  \newcommand{\frakt}[1]{\mbox{\Frakt #1}}      %
   \def\fr#1{\mathchoice{\fraktur {#1}}            
                        {\fraktur {#1}}            
                        {\fraktu {#1}}             
                        {\frakt {#1}}  }           

\newcommand{\Ss}{\fr S}

\def\db{{\mathchoice{\mbox{\rm db}}
                    {\mbox{\rm db}}
                    {\mbox{\scriptsize\rm db}}
                    {\mbox{\tiny\rm db}} }}

\def\id{{\mathchoice{\mbox{\rm id}}
                    {\mbox{\rm id}}
                    {\mbox{\scriptsize\rm id}}
                    {\mbox{\tiny\rm id}} }}

\def\part{\vdash}

\newcommand{\Phie}{{\Phi_e}}
\newcommand{\Phio}{{\Phi_o}}

\renewcommand{\Re}{{R_e}}
\newcommand{\Ro}{{R_o}}
\newcommand{\Go}{{G_o}}
\newcommand{\So}{{S_o}}
\newcommand{\Se}{{S_e}}
\newcommand{\Fe}{{F_e}}

\title{\mytitle\thanks{JOURNAL OF ALGEBRA to appear}}
\author{\myauthor\thanks{\mythanks}\\  \myaddress
\thanks{On leave from \myperaddress}\\ 
\myemail}

\date{}
\begin{document}
\maketitle
\begin{abstract}\myabstract\end{abstract}
\bibliographystyle{plain}

\section{Introduction}\label{sec1}

Quantum general linear super  groups were firstly studied by Manin. In \cite{manin2}, Manin introduced multiparameter deformations of the function ring M$(m|n)$ on the semi super group of super matrices. A quantum general linear super group is defined to consist of ``symmetries'' of a quantum space. Thus, Manin defined the quantum super semi-group M$_q(m|n)$ as a super bialgebra that  coacts universally on the quadratic algebra $A_q^{m|n}$ and its (graded) dual, which represent a quantum super space. M$_q(m|n)$ depends on $d(d-1)/2$ parameters $(d=m+n)$. By requiring that M$_q(m|n)$ should have the same Poincar\'e series as M$(m|n)$ does, Manin reduced the number of parameter to one continuous parameter $q$ and discrete parameters $\varepsilon_{ij}=\pm 1, \, 1\leq i, j\leq d $. Later, Demidov \cite{demidov} checked that $\varepsilon_{ij}$ should satisfy the condition $(\varepsilon_{ik}-\varepsilon_{ij})(\varepsilon_{jk}-\varepsilon_{ik})=0,\,  i<j<k$. The quantum super group GL$_q(m|n)$ is defined to be the Hopf envelope of M$_q(m|n)$. Sudbery had slightly changed Manin's construction and found a deformation depending on $n(n-1)/2+1$ continuous parameters \cite{sud1}. Sudbery's construction is also more compatible with the construction of Faddeev, Reshetikhin and Takhtajian \cite{frt}.
Multiparameter deformations for GL$(1|1)$ were also studied by several authors \cite{cor,dabr}.

The (super) determinant of a (super) matrix gives a criteria for the invertibility of the (super) matrix. For quantum super matrix it is not always possible to define the quantum super determinant. The standard way of obtaining a quantum (super) group from a given quantum (super) semi-group is to construct the Hopf envelope. Then one may use the quantum Koszul complex to define the quantum (super) determinant. The classical Koszul complex was introduced by Manin in \cite{manin2}. Lyubashenko \cite{lyu} and Gurevich \cite{gur} studied Koszul complexes associated to even (Hecke) symmetries. In \cite{ls}, Lyubashenko and Sudbery studied Koszul complexes associated to (non-even) closed Hecke symmetries. They gave an explicit formula of the quantum super determinant for the multiparameter deformation case. A formula for the quantum super determinant in this case was also given by Kobayashi and Uematsu without explicit explanations \cite{kob}.

It is still an open question if the Koszul complex always yields the quantum super determinant, i.e., if the complex has the cohomology group of dimension one, and if the quantum super determinant gives a criteria for the invertibility of the quantum super matrix.

In the classical case, the invertibility of a super matrix in a standard format $Z=\left( {A\atop C}{B\atop D}\right)$ is shown to be equivalent to the invertibility of the matrices $A$ and $D-CA^{-1}B$. The super determinant -- Berezinian -- is defined to be $\Ber Z=\det A\cdot\det (D-CA^{-1}B)$.
The aim of this paper is to show the q-analog of this fact, first for the multiparameter deformations and then for the deformations associated to Hecke sums of an odd and an even Hecke operator.

In $\S$\ref{sec2} and \S\ref{sec3} we recall the general construction of  quantum general linear super (semi) groups and the associated Koszul complex that induces the quantum super determinant. In $\S$\ref{sec6} we show that the quantum super determinant defined in \S\ref{sec3} has the same form as in the classical case and the quantum super matrix is invertible if this super determinant is invertible. We also prove the quantum analogs of the well-known formulas for the super determinant.  $\S$ \ref{sec4} is devoted to the study of deformations associated to Hecke sums of odd and even Hecke symmetries.

Throughout the paper, the parameter $q$ will be assumed to be invertible and non-root of unity other than 1.
 \section{Multiparameter Deformations of M$(m|n)$}\label{sec2}
Let $V$ be a finite-dimensional vector super space over $\kk$, a field of characteristic zero. Fix a homogeneous basis $x_1,x_2,\ldots, x_d$ of $V$, $d=\dim_\kk V$, and let $\hat{i}$ denote the parity (${\Bbb Z}_2$-gradation) of $x_i$. Let $\{p_{ij}, 1\leq i< j\leq d\}$ and $\{q_{ij}, 1\leq i< j\leq d\}$, be two sets of elements from $\kk^\times=\kk\setminus \{0\}$. For convenience, we shall also set $p_{ji}:=p_{ij}^{-1}, q_{ji}:=q_{ij}^{-1}, p_{ii}=q_{ii}=1$.
Let $R$ be an operator on $V\ot V$, whose matrix, with respect to the basis $x_i\ot x_j$, has the form:
\bba\nonumber &R(x_i\ot x_j)=x_k\ot x_lR^{kl}_{ij},&\\
\label{rmatrix} &R^{kl}_{ij}=\displaystyle
\frac{q-p_{ij}q_{ij}}{1+p_{ij}q_{ij}}\delta_{ij}^{kl}+(-
1)^{\hi\hj}\frac{p_{ij}(q+1)}{1+p_{ij}q_{ij}}\delta_{ij}^{lk},&\eea
here we use the convention of summing up by indices which appear in both upper and lower positions.

We shall assume that $p_{ij}q_{ij}$ satisfy the following conditions:
\bba\label{condition}
\begin{array}{l}
p_{ij}q_{ij}=q^{\va_{ij}}, \va_{ij}=\pm 1, i\neq j,\\
(\va_{ik}-\va_{ij})(\va_{jk}-\va_{ik})=0, \mbox{for } i<j<k.\end{array}\eea
One can easily verify that $R:V\ot V\lora V\ot V$ preserves parity, hence is a morphism in the category of super vector spaces.  \rref{rmatrix} implies that $R$ satisfies the Hecke equation $(R+1)(R-q)=0.$ The conditions in \rref{condition} ensure that $R$ satisfies the Yang-Baxter equation
\bba (R\ot\id_V)(\id_V\ot R)(R\ot \id_V)=(\id_V\ot R)(R\ot \id_V)(\id_V\ot R).\label{ybeq}\eea

 To $R$ there are associated the following quadratic algebras. The quadratic algebras $S$ and $\Lambd$ are factor algebras of the tensor algebra on $V$ by the following relations: 
\bba\begin{array}{ll}  \mbox{for  }S:& x_ix_jR_{kl}^{ij}=q x_kx_l\Longleftrightarrow x_ix_j=(-1)^{\hi\hj}q^{-1}_{ij}x_jx_i,\\
\mbox{for  }\Lambd :&x_ix_jR_{kl}^{ij}=-x_kx_l\Longleftrightarrow x_ix_j=-(-1)^{\hi\hj}p_{ij}x_jx_i.\end{array}\eea
 Let $E_1$ be a vector super space, generated over $\kk$ by a basis $\{z^i_j,i,j=1,2,\ldots,d\}$ with parities $ \hat{z}^i_j$ of $z_j^i$ equal to $\hat{i}+\hat{j}$.
The quadratic algebra $E$ is defined to be the factor algebra of the tensor algebra on $E_1$ by the ideal generated by the following relations:
\bba\label{zrelation} (-1)^{\hs(\hi+\hp)}R^{kl}_{ps} 
z^p_iz^s_j=(-1)^{\hl(\hq+\hk)}z^k_qz^l_nR_{ij}^{qn} .\eea
$E$ is a super bialgebra with the coproduct and the counit  defined as follows:
 \begin{equation}\label{multiplicative}
\Delta (z^j_i)=\sum_kz^j_k\ot z^k_i,\ \ \va (z^j_i)=\delta_i^j.\end{equation}
In fact, it can be shown that $E$ coacts universally on the super algebras $S$ and $\Lambd$ and the super bialgebra structure comes from the universality (see, e.g., \cite{manin2}, see also \cite{phdiss} for a more general setting). For convenience, we shall refer to a matrix, that satisfies the condition in \rref{zrelation}, as quantum super matrix, and, in the case all parameters are even, quantum matrix. A matrix of elements in $E$ is called multiplicative if it satisfies \rref{multiplicative}, in particular, a one-dimensional multiplicative matrix is a group-like element.

Assume that there are $m$ even parameters and $n$ odd parameters among $x_1,$ $ x_2,$ $ \ldots, x_d,$ $ m+n=d$, then $E$ is called the ``algebra of functions'' on M$_q(m|n)$ -- the multiparameter deformations  of the super semi-group M$(m|n)$.

 The algebras $S,\Lambd$ and $E$ were introduced in \cite{manin2}. $S$ and $\Lambd$ can be interpreted as the symmetric and exterior tensor algebras on the quantum super space while $E$ can be interpreted as the function algebra on the quantum  matrix super semi-group.

The matrix $R$ in \rref{rmatrix} is, so far, the most general $R$-matrix that induces a multiparameter deformation of GL$(m|n)$. It combines the deformation found by Manin \cite{manin2} and the deformation found by Sudbery \cite{sud1}. The algebra $E$ has a PBW-basis. In fact, one can show the following.

Assume that $R$ is given as in \rref{rmatrix}, then the following are equivalent:
\begin{itemize}
\item[1)] $R$ satisfies the Yang-Baxter equation.
\item[2)] $p_{ij},q_{ij}$ satisfy \rref{condition}.
\item[3)] $E$ has the correct Poincar\'e series and hence has a PBW-basis.
\end{itemize}

1) $\Longleftrightarrow$ 2) can be checked directly (cf. \cite{sud1}).

1) $\Longrightarrow$ 3) can be proved using the method of \cite{sud}, a more general result is obtained in \cite{ph97}. The existence of a PBW-basis can be shown using an algorithm of \cite{manin2}.

3) $\Longrightarrow$ 2) can be proved as in \cite{manin2}. The case $p_{ij}=q_{ij}$ has been checked by E. Demidov \cite{demidov}. 

Since $R$ satisfies the Yang-Baxter equation, it induces a coquasitriangular structure on $E$ and hence a braiding in the category of $E$-comodules. Moreover, $R$ is closed, i.e., there exists a matrix $S_{ij}^{kl}$ such that $R_{im}^{jn}S^{mk}_{nl}=S_{im}^{jn}R^{mk}_{nl}=\delta_{i}^{k}\delta^j_l\label{closed}$, hence one can extend the above category adding the dual of $V$. Recall that the dual super vector space to $V$ is a vector space $V^*$ generated by vectors $\xi^1,\xi^2,\ldots,\xi^a$ with parities $\hat{\xi}^i=\hat{i}$. The braiding on $V\ot V^*$ is given by means of $S$. For more details, the reader is referred to \cite{l-t,ls,gur}.

\hyphenation{Be-re-zin-ian}
\section{Quantum Super Group GL$_q(m|n)$ and Quantum Berezinian}\label{sec3}

Given a quantum semi (super) group, we can use Manin's construction to obtain the corresponding quantum (super) group. Manin's idea is to invert the multiplicative matrix of generators, in other words, to add to the set of generators the entries of the formal inverse matrix. Since we are working with non-commutative (super) algebras, the operation of taking inverse (i.e. the antipode) is not involutive. So that, in general, Manin's construction leads to an infinitely generated algebra. Thanks to the special property of the matrix $R$ (Yang-Baxter equation and closeness) the quantum super group is still finitely generated in our case. 

\subsection{The inverse of a quantum super matrix}\label{sec3.1}
 Let now $Z$ be a quantum super matrix in the sense of \S\ref{sec2}. Formally, we denote its inverse by $T$, thus we have relations like $t_k^jz^k_i=z_k^jt^k_i=\delta_i^j$, which imply that $T^t$ (but not $T$ !), up to some signs, is a quantum super matrix. 

The function ring $H$ on the quantum matrix super group GL$_q(m|n)$  is defined to be the Hopf envelope of $E$ \cite{manin2}. Let $\{t_i^j, i=1,2,\ldots, d\}$ be a set of homogeneous parameters of parities $\hat{t}^j_i=\hat{i}+\hat{j}$, $H$ is defined as a factor algebra of the free super algebra, generated by $z^j_i,t_i^j, i,j=1,2,\ldots, d$, by the ideal, generated by \rref{zrelation} and the following relations:
\bba\label{srelation} \sum_{k}(-1)^{\hi(\hk+1)}z^j_kt_i^k=\sum_{k}(-1)^{\hk(\hj+1)}t^j_kz_i^k=\delta_i^j,\eea
(cf. \cite{ls}).
Thus, we see that $(-1)^{\hk(\hj+1)}t^j_k$ is the inverse of $Z$. As in the non-super case (cf. \cite{frt,ph97b}), the antipode on $H$ can be given as follows \cite{phdiss}. Let $G$ be an $d\times d$ matrix, $G^j_i=S^{jl}_{il}$, then
\bba\label{antipodez}& S(z_i^j)=(-1)^{\hi(\hj+1)}t^j_i,&\\
 &\label{antipodet} S(t_i^j)=(-1)^{(k+i)j}{G_k^jz_l^kG^{-1}}_i^l.&
\eea
The matrix $G$ plays an important role in the theory and is called the quantum parity operator by Gurevich, we suppose the name reflection operator is more appropriate, since $G$ is the matrix for the canonical isomorphism $V\lora V^{**}$.  Direct computation shows that 
\bbas &G=\mbox{Diag} \{G_1,G_2,\ldots,G_{m+n}\},&\\
 &\Tr(G)=-[n-m]_q,&\eeas
where $G_k=-(-1)^{\hat{k}}q^{k_+-k_-}$, $k_+ (k_-)$ denotes the cardinality of even elements (odd elements) among $x_1,x_2,\ldots, x_k$.

The matrix $(t_i^j)^t$ is a multiplicative matrix: $\Delta t_i^j=\sum_kt_i^k\ot t_k^j.$ From \rref{srelation} and \rref{zrelation} we get
\bba\label{ztrelation}&
(-1)^{\hk(\hi+\hj)}R^{pj}_{ql}z_j^it_k^l=(-1)^{\hm(\hn+\hp)}t_n^pz_q^mR_{mk}^{ni},&\\
\label{ttrelation}&
(-1)^{\hs(\hi+\hp)}R^{kl}_{ps} 
t^s_jt^p_i=(-1)^{\hl(\hq+\hk)}t^l_nt^k_q R_{ij}^{qn}.&\eea

Remark that the construction of this section applies well for any closed Yang-Baxter operator. In \S4 we shall study deformations associated to the Hecke sum of an odd and an even Hecke operator.

\subsection{The Quantum Determinant}\label{sec3.2}
In the non-super case, the quantum determinant is used to construct a Hopf algebra from the matrix bialgebra. It is a group-like element in the bialgebra $E$, which has a special property, so that we can localize $E$ with respect to this element to get a Hopf algebra. If it commutes with all elements we can set it equal to 1 to get a deformation of the special linear group. Quantum determinant can be defined for any even Hecke operators \cite{gur}.

Since, later, we will want to express the quantum super determinant as a product of quantum determinant of certain quantum matrices, we shall need to know, for what kind of matrix of non-commuting entries it is possible to define a quantum determinant.

 Let us recall the quantum determinant for the multiparameter deformation of GL$(d)$. Thus assume that $\hi =0 $ for all $i$, in other words $n=0, m=d$. Then $\Lambd_d$ (the $d^{th}$ homogeneous component of $\Lambda$) is one-dimensional and can be generated by  $x_1x_2\cdots x_d$. Since $\Lambd_d$ is an $E$-comodule, the coaction $\delta 
:\Lambd_d\lora\Lambd_d\ot E$ defines a group-like element in $E$:
$ \delta(x_{1}x_{2}\cdots x_{d})=x_{1}x_{2}\cdots x_{d}\ot \det_qZ.$
Let us denote for any element $\si$ of the symmetric group $\Ss_d$, $\si(p):=\prod_{i<j,\ i\si>j\si}p_{i\si,j\si}.$
Then we have $\det_qZ=\sum_{\si\in\Ss_d}\si(-p)z_{1}^{1\si}z_{2}^{2\si}\cdots z_{d}^{d\si}$. More generally, since any element $x_{1\tau}x_{2\tau}\cdots x_{d\tau}$, $\tau\in\Ss_d$, generates $\Lambda^d$ as well, we have
\bba\label{dete}\det_qZ=\sum_{\si\in\Ss_d}\textstyle\frac{\si(-p)}{\tau(-p)}z_{1\tau}^{1\si}z_{2\tau}^{2\si}\cdots z_{d\tau}^{d\si},\eea
for all $\tau$ from $\Ss_d$.
Let $H:=E[\det_qZ^{-1}]$ then $H$ is a Hopf algebra \cite{gur,frt}. Let $S$ be the antipode and $t_i^j:=S(z_i^j)$, then, by definition of the antipode, 
\begin{equation}\label{eq5a}
\sum_k z^j_kt^k_i=\sum_k t^j_kz^k_i=\delta_i^j.\end{equation}
\rref{eq5a} and \rref{zrelation} imply the commuting rule between $z_i^j$ and $t_k^l$:
\bbas R^{pj}_{ql}z_j^it_k^l=t_n^pz_q^mR_{mk}^{ni},\\
\label{ttrel}R_{sp}^{lk} t_k^it_l^j=t_p^mt_s^nR^{ji}_{nm}.\eeas
Let $\Lambda^\vee$ be the dual quadratic algebra to $\Lambda$ (see the next section). Then $\Lambda^\vee$ is also an $H$-comodule. $\Lambda^\vee_d$ is of dimension one, and hence defines a group-like element in $H$, which is the quantum determinant of $T=(t_i^j)$. Notice that $\Lambda^\vee_d$ is dual to $\Lambda_d$ as $H$-comodules, thus
\bba\label{invdet} \det_q Z=(\det_qZ^{-1})^{-1}.\eea

If all parameters are odd, i.e., $\hi =1$ for all $i$, then $S_d$ will be 
one-dimensional and generated by any of the elements $x_{1\tau}x_{2\tau}\cdots x_{d\tau}$. Denote by 
$\det_qZ$ the quantum determinant: $\delta (x_{1\tau}x_{2\tau}\cdots x_{d\tau})=x_{1\tau}x_{2\tau}\cdots x_{d\tau}\otimes\det_qZ$. Thus, we have
\bba\label{deto}
\det_qZ=\sum_{\si\in\Ss_d}\textstyle\frac{\si(q^{-1})}{\tau(q^{-1})}z_{1\tau}^{1\si}z_{2\tau}^{2\si}\cdots z_{d\tau}^{d\si}.
\eea
And hence, for $H:=E[\det_qZ^{-1}]$ we have the analogous assertions.

 Assume again that all elements $x_i$ are even. Let $K$ be an algebra and $(a_i^j)_{i,j=1}^{d}$ be a matrix of elements from $K$. $(a_i^j)$ is called a quantum matrix if there exists an algebra homomorphism, denoted by $A$, $A:E\lora K, z^j_i\loma a_i^j$, in other words, $A$ is a $K$-point of $\mbox{M}_q(d)$. Under $A$, the quantum determinant
 $\det_q$ is transposed into a certain element in $K$, denoted by $\det_qA$. From the discussion above we immediately have the  following lemma, which will be of use in the next sections. 
\begin{lem}\label{trans} Let K be an algebra and $(a_i^j)_{i,j=1}^{d}$ be a quantum matrix in $K$. Then there exist elements $\det_qA$ and $b_i^j$ in $K$, such that
\bbas \sum_ka_k^jb_i^k=\sum_kb_k^ja_i^k=\det_qA.\eeas
If moreover  $\det_qA$ is invertible, then there are elements $c_i^j$ in $K$ such that
\bbas \sum_ka_k^jc_i^k=\sum_kc_k^ja_i^k=\delta_i^j,\eeas
and $c_i^j$ obey the relation:
\bbs c_k^ic_l^jR_{sp}^{lk}=R^{ji}_{nm}c_p^mc_s^n.\ees
Hence $\det_q C$ can be defined and $\det_q C=(\det_q A)\iii$, in other words, 
$\det_q (A\iii)=(\det_q A)\iii.$
\end{lem}

\subsection{Koszul complex and Berezinian}\label{sec3.3}
Generally, finding a group-like element of a bialgebra is equivalent to finding a one-dimensional comodule. Unlike in the non-super case, in the super case, the algebras $S$ and $\Lambd$ are both infinite-dimensional, hence we cannot derive from them group like elements as in the previous section. The quantum Berezinian was introduced by Manin \cite{manin2} as the homology of certain complex defined on $S$ and $\Lambd$. It is not contained in $E$ but rather in $H$. The explicit formula for the quantum determinant was given by Lyubashenko and Sudbery \cite{ls}.

By definition $V$ is an $H$-comodule with the coaction $\delta :V\lora V\ot H 
\delta (x_i)=x_j\ot z_i^j$. The dual space $V^*$ is also  an $H$-comodule, with the coaction $\delta (\xi^i)=\xi^j\ot t_j^i$. The algebras $\Lambd$ and $S$ have been defined by
\bbas \Lambd=T(V)/<\Im(R+1)>,\ \ \ S=T(V)/<\Im(R-q)>.\eeas
Let us consider the following pairing on $V^*\ot V^*$ and $V\ot V: <\xi^i\ot\xi^j|x_k\ot x_l>:=\delta_l^i\delta^j_k.\label{pairing}$
It makes $V^*\ot V^*$  a dual object to $V\ot V$ in the category of super vector spaces.
Let $R^*$ be the dual operator to $R$ with respect to this pairing. We have $ R^*(\xi^i\ot \xi^j)=\xi^k\ot \xi^lR^{ji}_{lk}.$
 Thus we can define quadratic algebras on $V^*$: $S^\vee$ and $\Lambd^\vee$ by:
\bbas \Lambd ^\vee:=T(V^*)/<\Im(R^*+1 )>,\ \ \
S^\vee:=T(V^*)/<\Im(R^*-q)>.\eeas
 The relations on $\Lambd^\vee$ and $S^\vee$ are 
$ \xi^k\xi^l=-(-1)^{\hk\hl}q^{-1}_{kl}\xi^l\xi^k$ and $\xi^k\xi^l=(-1)^{\hk\hl}p_{kl}\xi^l\xi^k$, respectively.

Set ${\bf K}^{k,l}:=\Lambd_k\ot S^\vee_l$. Let $\d:{\bf K}^{k,l}\lora {\bf K}^{k+1,l+1} $ be  the morphism
 \bbas\label{differential} \d=\id_V^{\ot k}\ot \db\ot \id_{V^*}^{\ot l}:V^{\ot k}\ot V^{*\ot l}\lora V^{\ot k+1}\ot V^{*\ot l+1},\eeas
where $\db$ is the morphism ${\kk}\lora V\ot V^*$, $\db(1_{\kk})=\sum_ix_i\ot\xi^i.$
Then $\d^2=0$. Hence for every $k$, $\d$ induces a differential in the complex
\bbs{\bf K}_k^\bullet=\cdots \stackrel{\d}{\lora} {\bf K}^{0,k}\stackrel{\d}{\lora} {\bf K}^{1,k+1}\stackrel{\d}{\lora}\cdots\ees
Therefore ${\bf K}_k^{\bullet}$ 
is a complex in the category of $H$-comodules,  its cohomology group is an $H$-comodule. Assume that there are among $x_1,x_2,\ldots, x_d,$ $m$ even elements and $n$ odd elements, $m+n=d$. Reordering the $x_i$'s we can assume that $\hat{i}=0$ for $i\leq m$ and $\hat{i}=1$ for $i> m$.
\begin{edl}\cite{manin2}\label{manin} 
The complex $({\bf K}_k^{\bullet},\d)$ is exact everywhere  
except at the term ${\bf K}_{n-m}^{m,n}$, and the cohomology group of ${\bf K}_{n-m}^\bullet$ is one-dimensional, generated by the element
\bba\label{cohomo}x_1x_2\cdots x_m\xi^{m+1}\xi^{m+2}\cdots \xi^{m+n} \mbox{  \nml  modulo }\Im \d.\eea
\end{edl}
The group-like element in $H$ induced from this cohomology group is called the quantum super determinant or quantum Berezinian.
\begin{edl}\cite{ls} The quantum Berezinian can be given by the following formula, which is valid for any $\tau\in\Ss_m$ and $\theta\in\Ss_n$, 
\bba\label{sdet}
{\sdet Z} 
\displaystyle=\sum_{\si\in\Ss_m}\textstyle\frac{\si(-p)}{\tau(-p)}
z_{1\tau}^{1\si}z_{2\tau}^{2\si}\cdots z_{m\tau}^{m\si}\times
 \displaystyle\sum_{\nu\in\Ss_n}\textstyle\frac{\theta(-p,m)}{\nu(-p,m)}
t^{m+1\theta}_{m+1\nu}t^{m+2\theta}_{m+2\nu}\cdots t^{m+n\theta}_{m+n\nu},\eea
where $\nu(p,m):=\prod_{i<j,i\nu> i\nu}p_{m+i\nu,m+j\nu}$.
\end{edl}
The  complex in our form here was given by Lyubashenko and Gurevich \cite{lyu,gur}, it slightly differs from one introduced by Manin \cite{manin2} (since the pairings were defined differently).

\section{Structure of Quantum Super Matrix}\label{sec6}
In \S\ref{sec2} we have  defined the structure of the function ring $H$ on the quantum group GL$_q(m|n)$ as well as the quantum Berezinian. However, the Hopf algebra $H$ is still generated by too many generators ($2d^2$). The aim of this section is to reduce the number of generators to minimal ($d^2+2$) as well as provide quantum analogs of the well-known formulas for the classical Berezinian.
 
 We will assume that $\daysub{x}{m}$  are even parameters and $x_{m+1},\cdots 
,x_{m+n}$ are odd parameters. The matrices $Z=(z_i^j)_{i,j=1}^{m+n}$ and 
$T=(t_j^i)_{i,j=1}^{m+n}$ decompose into 4 matrices:
\bbas Z=\left(\begin{array}{ll}A & B\\ C&D\end{array}\right),\ \ \ 
T=\left(\begin{array}{ll}X & Y\\ U&V\end{array}\right) .\eeas
The antipode is then:
\bbas S\left(\begin{array}{lr}A & B\\ C&D\end{array}\right)=
\left(\begin{array}{lr}X & -Y\\ U&V\end{array}\right) .\eeas
Thus we have:
\bba\label{santipode1}\left(\begin{array}{ll}A & B\\ C&D\end{array}\right)
\left(\begin{array}{lr}X & -Y\\ U&V\end{array}\right)=\left(\begin{array}{ll}1 & 
0\\ 0&1\end{array}\right) ,
\left(\begin{array}{ll}X & -Y\\ U&V\end{array}\right)
\left(\begin{array}{lr}A & B\\ C&D\end{array}\right)=\left(\begin{array}{ll}1 & 0\\ 0&1\end{array}\right)\eea
According to \rref{sdet}:
\bbas \sdet Z=\det_qA\cdot \det_qV.\eeas

According to \rref{rmatrix} and \rref{ztrelation}, we have the following commuting rule between $z_i^j$ and $t_k^l$, $i,j=\dayso{m}$, $k,l=m+1,m+2,\cdots ,m+n$:
\bbas 
\frac{p_{il}}{1+p_{il}q_{il}}z_i^jt_k^l=\frac{p_{jk}}{1+p_{jk}q_{jk}}t_k^lz_i^j,
\eeas
from which it follows immediately that $\det_qA$ and $\det_qV$ commute.
In fact, every term of $\det_qA$ in \rref{dete} commutes with every term of $\det_qV$ in \rref{deto}.
Since $\sdet Z$ is a group-like element, it is invertible in $H$. Therefore 
$\det_qA$ and $\det_qV$  are invertible. According to Lemma \ref{trans}, we can define matrices $\Ai$ and $\Vi$, such that $\Ai A=A\Ai =1,\ \Vi V=V\Vi=1$.

Multiplying both sides of the first equation in \rref{santipode1} by Diag$\{\Ai ,1\}$ from the left and by 
Diag$\{1,\Vi\}$ from the right, we obtain
\bba\label{equation}\left(\begin{array}{ll}X+\Ai BU & -Y\Vi+\Ai B\\ 
CX+DU&-CY\Vi+D\end{array}\right)=
\left(\begin{array}{ll}\Ai & 0\\ 0&\Vi\end{array}\right).\eea
Hence $H:=D-C\Ai B=\Vi$. Therefore the quantum Berezinian is
\bba\label{qtber1} \Ber_q Z=\det_qA\cdot \det_q(D-C\Ai B)^{-1}.\eea
According to \rref{ttrelation}, $V$ is a quantum matrix. Using Lemma \ref{trans} we obtain the  relation \rref{zrelation} for the elements $h_i^j$ of $H$:
\bbas{R_o}^{ps}_{kl} h^k_ih^l_j=h^p_qh^s_n{R_o}_{ij}^{qn}, \eeas
where $R_o$ denotes the sub-matrix of $R$ of  elements, whose parities of both 
indices are odd. Thus $H$ is a quantum matrix. Let $\det_qH$ be its quantum determinant.
From \rref{santipode1} and \rref{equation} we obtain $ Y=\Ai B\Hi ,  U=-H\iii  C\Ai.$
Whence $ X=\Ai +\Ai B\Hi C\Ai.$ Thus
\bba\left(\begin{array}{lr}X&Y\\ U&V\end{array}\right)=
\left(\begin{array}{lr}\Ai +\Ai B\Hi C\Ai&  \Ai B\Hi\\    -\Hi C\Ai& 
\Hi\end{array}\right).\eea

Let us consider another Koszul complex ${\bf L}^{\bullet\bullet}:=S\ot \Lambd^\vee$, the cohomology 
of which is also one-dimensional, hence induces a group-like element
\bbas \det_qD\cdot\det_qX=\det_qX\cdot\det_qD.\eeas
\begin{lem}\label{lemisdet} The following equations hold
\bbas &S(\det_qA)=\det_qX,& \\
 &S(\det_qV)=\det_qD.&\eeas\end{lem}
\proof Since $S(A)=X$, let $w=w_m$ be the permutation that reverses the order of $1,2,\ldots,m$, we then have
\bbs \frac{w\si(-p)}{\si(-p)}=\prod_{1\leq i<j\leq m}(-p_{ij}),\quad \forall \si\in\Ss_m.\ees
Therefore, by \rref{dete}, we have 
\bbas S(\det_qA)&=& S\left( \sum_{\si\in\Ss_m}\textstyle\frac{\si(-p)}{\tau(-p)}z_{1\tau}^{1\si}z_{2\tau}^{2\si}\cdots z_{m\tau}^{m\si}\right)\\
 &=&\sum_{\si\in\Ss_m}\textstyle\frac{\si(-p)}{\tau(-p)}t_{m\tau}^{m\si}\cdots t_{2\tau}^{2\si} t_{1\tau}^{1\si}\\
 & =& \sum_{\si\in\Ss_m}\textstyle\frac{w\si(-p)}{w\tau(-p)}t_{1w\tau}^{1w\si}t_{2w\tau}^{2w\si}\cdots z_{mw\tau}^{mw\si}\\  &=& \det_qX.\eeas

According to \rref{antipodet}, $S(V)=GDG^{-1}$, where $G$ is a diagonal matrix. Therefore $\det_q(GDG^{-1})=\det_qD$. Hence, $S(\det_qV)=\det_q(GDG^{-1})=\det_qD$.\eee

\begin{pro}\label{isdet} $(\sdet Z)^{-1}=\det_qD\cdot 
\det_q(A-BD^{-1}C)^{-1}.$\end{pro}
\proof Since $\Ber_qZ$ is a group-like element, $S(\Ber_q Z)= (\Ber_q Z)^{-1}$. On the other hand, we have $ X=(A-B\Di C)^{-1}.$\eee
The formula in Proposition \ref{isdet} gives us another formula for the quantum Berezinian:
\bbas \Ber Z=\det_q\Di\det_q(A-B\Di C) .\eeas

We summarize the results obtained in a theorem.
\begin{edl}\label{firstimp} Let $Z$ be a quantum super matrix defined as in \rref{zrelation}, assume that $Z$ is written in the standard format:
\bbs
Z=\left(\begin{array}{ll}A & B\\ C&D\end{array}\right) .\ees
 Then $A$ and  $D$ are quantum matrices. If $Z$ is invertible, then $A$ and $D$ are invertible too. Moreover, the matrices $A-BD^{-1}C$ and $D-CA^{-1}B$ are invertible, too, and the quantum Berezinian of $Z$ can be given by
\bbs
\Ber_q Z=\det_qA\cdot \det_q(D-C\Ai B)^{-1}=\det_q\Di\det_q(A-B\Di C).\eqeee
\ees
\end{edl}

It is now natural to ask, if the converse assertion to Theorem \ref{firstimp} holds. Below we will show that the answer is positive.

If $\varepsilon_{ij}$ satisfy the condition in \rref{condition} then we can reorder $x_i$ in such a way that $\varepsilon_{ij}=1$ for $i<j$. However, in this case we may not have the standard form for $Z$. We shall prove here the converse assertion to Theorem (2.6), provided that $\varepsilon_{ij}=1$ for $i<j$ and $Z$ is in the standard format. The proof for the general case can be done in the same way. However, in the general case, the relations become more complicated.

If we rearrange the indices of $R=R_{ij}^{kl}$ in a form that is compatible with the decomposition of $Z$ into four matrices, we shall have the following form of $R$:

\begin{equation}\label{eq0}
R = \left| \begin{array}{cccc} 
R_e & 0 & 0 & 0\\
0 & 0 & Q & 0\\
0 & P & (q-1)I & 0\\
0 & 0 & 0 & R_o
\end{array} \right| \end{equation}
where $R_e=(R)_{i,j=1}^m, R_o=(R)_{i,j=m+1}^{m+n},\displaystyle Q^{kl}_{ji}= q_{ij}\delta_{ij}^{kl}, P^{lk}_{ij}=p_{ij}\delta_{ij}^{kl}, i,k\leq m<j,l,$ whence $ PQ=qI$. 

For convenience, we shall denote $T_1=T\otimes I, T_2=I\otimes T$ where $T$ is a (not necessarily square) matrix, usually of the form $m\times n, m\times n, n\times m $ or $n\times n$, $I$ is the identity matrix of the form $m\times m$ or $n\times n$. The choice of $I$ will be taken appropriately in each context. More generally, we shall also consider the matrix $T_3=I\otimes I\otimes T$, and so on. 

The relations in (\ref{zrelation}) can be rewritten in the following form:
\begin{equation}\label{eq1}\begin{array}{rclcrcl}
R_e A_1 A_2 & =& A_1A_2R_e &  & R_eB_1B_2 & =&  B_1B_2R_o\\
C_1C_2R_e &=&  R_oC_1C_2 & &  R_oD_1D_2  &= & D_1D_2R_o\\
R_eA_1B_2 & =& B_1A_2P &  & A_1C_2R_e  &= &qP^{-1}C_1A_2\\
 R_oC_1D_2 & =&  -D_1C_2P&  & -B_1D_2R_o & =&  qP^{-1}D_1B_2 \\
\multicolumn{3}{r}{ qA_1D_2P^{-1}- qP^{-1}D_1A_2}&  =& \multicolumn{3}{l}{ (q-1)B_1C_2 }\\
\multicolumn{3}{r}{ B_1C_2P } &= &\multicolumn{3}{l}{ -qP^{-1}C_1B_2 .}\end{array}
\end{equation} 
\begin{lem}\label{lem1}
Assume that $A$ is invertible. Let $H=D-CA^{-1}B$ then $H$ is a quantum matrix: $R_oH_1H_2=H_1H_2R_o$.
\end{lem}
\proof
First we have 
\[ \begin{array}{rcl}
(CA^{-1}B)_1(CA^{-1}B)_2 &=& -qC_1A_1^{-1}P^{-1}C_1B_2P^{-1}A_2^{-1}B_2 \\
& =& -C_1C_2R_eA_2^{-1}A_1^{-1}R_e^{-1}B_1B_2\\
& = & C_1C_2A_2^{-1}A_1^{-1}B_1B_2. \end{array} \]
Therefore
\[ \begin{array}{rcl}
R_oC_1A_1^{-1}B_1C_2A_2^{-1}B_2 &=& R_oC_1C_2A_2^{-1}A_1^{-1}B_1B_2\\
& =& C_1C_2A_2^{-1}A_1^{-1}B_1B_2R_o\\
& =& C_1A_1^{-1}B_1C_2A_2^{-1}B_2R_o.\\
\end{array} \]
On the other hand, we have 
\[ \begin{array}{rcl}
D_1C_2A_2^{-1}B_2 &=& -R_oC_1D_2P^{-1}A_2^{-1}B_2\\
&=& -R_oC_1(A_1^{-1}P^{-1}D_1+(q-1)q^{-1}A_1^{-1}B_1C_2A_2^{-1})B_2\\
& =& q^{-1}R_oC_1A_1^{-1}B_1D_2R_o-(1-q^{-1})R_oC_1A_1^{-1}B_1C_2A_2^{-1}B_2.
\end{array} \]
The second term, as shown above, commutes with $R_o$. Thus, we have
\begin{eqnarray*}
\lefteqn{R_o(D_1C_2A_2^{-1}B_2+C_1A_1^{-1}B_1D_2) }\\
 &= & q^{-1}R_o^2C_1A_1^{-1}B_1D_2R_o+R_oC_1A_1^{-1}B_1D_2-(1-q^{-1})R_o^2C_1A_1^{-1}B_1C_2A_2^{-1}B_2\\
& = &C_1A_1^{-1}B_1D_2R_o+R_oC_1A_1^{-1}B_1D_2(1+(1-q^{-1})R_o)\\
 &&-(1-q^{-1})R_o^2C_1A_1^{-1}B_1C_2A_2^{-1}B_2\\
 &= & (D_1C_2A_2^{-1}B_2+C_1A_1^{-1}B_1D_2)R_o
\end{eqnarray*}
Since $D$ also satisfies $R_0D_1D_2=D_1D_2R_o$, the lemma is proved.\eee

\begin{cor}\label{cor2}
 If $A$ is invertible then we can define the quantum Berezinian for $Z$:
\[ Ber Z=\det_q A^{-1}\det_q (D-CA^{-1}B).\eqeee\]\end{cor}
The analogous assertion to Corollary \ref{cor2} holds if we assume $D$ to be invertible. If we assume that both $A$ and $D$ are invertible we will have an element $\det_q D^{-1}\det_q (D-CA^{-1}B)$ which is expected to be the inverse of Ber $Z$ (cf. Proposition \ref{isdet}). Since we don't have the Hopf algebra structure yet, the method of the previous section cannot be applied. A direct verification seems to be very complicated. We prove instead that
 if $A$ and $D$ are invertible then $Z$ is invertible.

According to Corollary \ref{cor2} it is sufficient to show that $D-CA^{-1}B$ is invertible if $A$ and $D$ are, or equivalently, to show that $CA^{-1}BD^{-1}$ is nilpotent. This fact follows from the lemma and the corollary below.
 
\begin{lem}\label{lem4}
The matrix $K=CA^{-1}BD^{-1}$ satisfies the following relation 
\[ K_1(R_oT)_1K_2(R_oT)_2\cdots (R_oT)_{mn}K_{mn+1}=0,\]
 where $T$ is the usual permuting operator on $V_{\overline{0}}\otimes V_{\overline{1}}: T_{ij}^{lk}= \delta_{ij}^{kl}, i,j\leq m<k,l$.
\end{lem}
\proof According to \rref{eq1},
\bbas &D_1^{-1}R_oTC_2=-C_2PTD^{-1}_1 \quad qA^{-1}P^{-1}T^{-1}C_2=C_2\Re T^{-1}A^{-1}_1& \\
& -B_1C_2PT=qP^{-1}T^{-1}C_2B_1.& \eeas
Therefore
\begin{eqnarray*}
C_1A_1^{-1}B_1D_1^{-1}R_oTC_2 &=& -C_1A_1^{-1}B_1C_2PTD_1^{-1}\\
&=& qC_1A_1^{-1}P^{-1}T^{-1}C_2B_1D_1^{-1}\\
& =& C_1C_2R_e^{-1}A_1^{-1}B_1D_1^{-1}. \end{eqnarray*}
Analogously, we can show that
\[ K_1(R_oT)_1K_2(R_oT)_2\cdots (R_oT)_{k-1}K_{k}\]
can be led to the form that contains $C_1C_2\cdots C_k$. Since  $C_1C_2\cdots C_{mn+1}=0$, the assertion is proved.\eee
\begin{cor}\label{cor5}$K_1K_2\cdots K_{mn+1}=0$.\end{cor}
\proof We prove, that if $K_1(R_oT)_1K(R_oT)_2\cdots (R_oT)_{l-1}K_{l}=0$, then $K_1K_2\cdots K_l$ $= 0$. We shall prove this for $l=2$, the general case being similarly proved. Thus, assume $K_1(R_oT)_1K_2=0,$ i.e.
\[\sum_{m+1\leq m, q\leq m+n} K_m^iK_l^q \: \Ro_{qk}^{mj}=0,\quad \forall m+1\leq i,j,k,l\leq m+n.\]
Since the matrix $\Ro$ is closed, i.e., there exists a matrix $\So$, such that $\Ro^{im}_{kn}\So_{ml}^{nj}=\So^{im}_{kn}\Ro_{ml}^{nj}=\delta^i_l\delta^j_k$, from the above equation we derive
\bbas K^i_jK^k_l=0, \quad \forall m+1\leq i,j,k,l\leq m+n.\eqeee\eeas

 Corollary \ref{cor5} shows that $K$ is a nilpotent matrix. Hence $I-K=I-CA^{-1}BD^{-1}$ is invertible. Analogously we have $I-D^{-1}CA^{-1}B$ is invertible. Therefore $D-CA^{-1}B$ is invertible. Thus we have proved
\begin{edl}\label{edl6} 
Let $R$ be the matrix in \rref{rmatrix} with $\varepsilon_{ij}=1, \forall i<j,$ then the function algebra on the associated quantum general linear super group GL$_R(m|n)$ can be characterized as:

1. the algebra generated by entries of $A, B, C, D$ and $\det_q A^{-1}, \det_q(D-CA^{-1}B)^{-1}$, factorized by the relations in (\ref{eq1}).

2. the algebra generated by entries of $A, B, C, D$ and $\det_qA^{-1},\det_q D^{-1}$, factorized by the relations in (\ref{eq1}).
\end{edl}

\noindent{\bf A remark on the general case.} In the general case, the matrix $R$ has the following form:
\bbs
R = \left| \begin{array}{cccc} 
R_e & 0 & 0 & 0\\
0 & R_1 & Q & 0\\
0 & P & R_2 & 0\\
0 & 0 & 0 & R_o
\end{array} \right| \ees
where the matrices $R_1$ and $R_2$ satisfy certain conditions: $(R_1Q)^2=(q-1)R_1Q,$ $R_1Q+QR_2=(q-1)Q$. The relations in \rref{eq1} become more complicated. To prove Lemma \ref{lem1} one should first establish the relations between $H$ and $A,B,C,D$, which can be predicted from \rref{ztrelation}. More work is needed for Lemma \ref{lem4}.

\section{Generalization}\label{sec4}
In this section we generalize the results obtained in the previous sections for the Hecke sum of an even and an odd Hecke operator. The decomposition of $R$ in \rref{eq1} is a special case of the Hecke sum, first introduced by Gurevich \cite{gur}, and studied by Lyubashenko and Sudbery in \cite{ls}.

By definition, an even Hecke operator or an even Hecke symmetry on a non-super vector space $V$ is an operator on $V\otimes V$, satisfying the Yang-Baxter and Hecke equations, being closed (see page \pageref{ybeq}) and such that the associated anti-symmetric tensor algebra $\Lambda$ is finite dimensional. The degree of the highest non-vanishing homogeneous component of $\Lambda$ is called  the rank of $R$. It was shown by Gurevich that for such an operator, the associated bialgebra, defined by relations in (\ref{zrelation}), possesses a group-like element that plays the role of the quantum determinant. The matrix quantum group associated to an even Hecke operator was studied in \cite{gur, ph97b}. An operator $R$ on $V\otimes V$ is an odd Hecke operator if $-qR^{-1}$ is an even Hecke operator. The rank of $R$ is defined to be the rank of $-qR^{-1}$.   

Let $R_e$ and $R_o$ be an even and an odd Hecke operator on $V_0$ and $V_1$ respectively. Let $P$ be an invertible, closed operator $P: V_1\otimes V_0\lora V_0\ot V_1$ such that
\begin{eqnarray}\nonumber
 &P_2P_1\Re_{2} = \Re_1P_2P_1,&\\
 &P_1P_2\Ro_1  =  \Ro_{2}P_1P_2.&\label{eqforp} \end{eqnarray}
Let $V$ be the super vector space, such that $V_{\overline{0}}=V_0, V_{\overline{1}}=V_1$ and let $R$ be an operator on $V\otimes V$ defined from $R_e, R_o$ and $P$ as in \rref{eq0}, where $Q=qP^{-1}$. Then $R$ is a closed Hecke operator. Along the same line as in \S\ref{sec2} and \S\ref{sec3} we can define the super bialgebra $E$ and its Hopf envelope $H$ associated to $R$, which are called the matrix quantum super (semi) group of type $(r, s)$, where $r$ (resp. $s$) is the rank of $R_e$ (resp. $R_o$).

The following questions are to be answered:
\begin{itemize}
\item[1)] What is the Poincar\'e series of $E$ and when does $E$ have a PBW-basis?
\item[2)] Can one define the quantum (super) determinant?
\item[3)] Does Eq. \rref{qtber1}   hold?
\item[4)] Is Theorem \ref{edl6}  still true?
\end{itemize}

The Poincar\'e series of $E$ can be given in terms of the Poincar\'e series of the bialgebras $E_0$ defined on $V_0$ by $\Re$ and $E_1$ defined on $V_1$ by $R_o$ \cite{ph97c}
\bba\label{poincare}
P_E(t)=P_{E_0}(t)*P_{E_1}(t).\eea

\begin{lem}\label{lem41}
The following equations hold: \[ \begin{array}{c} B_1B_2\cdots B_{rs+1}=0,\\
C_1C_2\cdots C_{rs+1} = 0.\end{array}\]\end{lem}
\proof
We give the proof of the first equation, the proof for the second one is similar. Consider the algebra generated by $R_{ei}, i=1,\cdots, rs+1$ (resp. $R_{oi}$), which is a factor of the Hecke algebra $\H_{rs+1}$. $\H_{rs+1}$ is semi-simple, its central primitive idempotents are indexed by $\{ \lambda |\lambda \vdash rs+1\}$ \cite{dj1}. Let $R_{e\lambda}$ (resp. $R_{o\lambda}$) be the image of the central primitive idempotent in $\H_n$ corresponding to $\lambda$ in the algebra generated by $R_{ei}$ (resp. $R_{oi}$). Thus, we have $1=\sum_{\lambda\vdash rs+1}R_{e\lambda}$ and the analogous equation for $R_{o\lambda}$. 

The relation in \rref{eq1} for $B$ implies:
\[ R_{e\lambda}B_1B_2\cdots B_{rs+1}=B_1B_2\cdots B_{rs+1} R_{o\lambda},
\]
  hence
\[
B_1B_2\cdots B_{rs+1}=\sum_{\lambda \vdash rs+1}R_{e\lambda}B_1B_2\cdots B_{rs+1}=\sum_{\lambda\vdash rs+1}B_1B_2\cdots B_{rs+1} R_{o\lambda}.
\]
 Since $R_e$ is an even Hecke operator of rank $r, e_\lambda=0$ unless $\lambda_i=0, \forall i\geq r+1$, analogously, $o_\lambda=0$ unless $\lambda_1\leq s$ \cite{ph97b}. Whence for all $\lambda\vdash (rs+1),$ either $\Re_\lambda=0$ or $\Ro_\lambda=0$. The equation above then implies that $B_1B_2\cdots B_{rs+1}=0$.\eee

The relations in \rref{eq1} for the matrices $A, B, C, D$ 
allow us to express any element of $E$ as a linear combination of elements of the form $d c b a$, where $d $  (resp. $ c, b, a)$ is the monomial in the entries of $D$ (resp. $C, B, A$). Hence, if $\{a_i|i\in I\}$ (resp. $\{d_j|j\in J\}$) form a basis of the subalgebra of $E$ generated by entries of $A$ (resp. $D$), i.e. $E_0$ (resp. $E_1$), then $d_jc_kb_la_i$ span $E$, where $c_k$ (resp. $b_l$) is the ordered monomials in the entries of $C(D)$. From \rref{poincare} and Lemma \ref{lem41} it follows that $d_jc_kb_la_i$ form a basis for $E$. Thus we proved
\begin{pro}\label{pro2}
If on $E_0$ and $E_1$ there exist PBW-bases then there exists such a basis on $E$.\end{pro}

Lyubashenko and Sudbery \cite{ls} showed that the cohomology group of the Koszul complex associated to $R$ is isomorphic to the tensor product of the cohomology groups of the Koszul complexes associated to $R_e$ and $R_o$. Since $R_e$ is an even Hecke operator, the cohomology group for $R_e$ is isomorphic to $\Lambda^r_0$, where  $\Lambda_0$  is the quantum antisymmetric tensor algebra defined on $V_0$, $r$ is the rank of $R_e$. Analogously, the cohomology group for $R_o$ is isomorphic to $S_1^{\vee s}$, where $s$ is the rank of $R_o$. Hence, the cohomology of the Koszul complex on $V$ is isomorphic to $\Lambda_0^r\otimes S_1^{\vee s}$. This is a generalization of Theorem \ref{manin}. Thus, we can define the quantum Berezinian in the Hopf envelope $H$ of $E$. 
\begin{pro}\label{pro43}
 Let $\det_qA$ (resp. $\det_qV)$ be the quantum determinant of $A$ (resp. $V$), then they commute and the quantum Berezinian of $Z$ is given by 
\[ \Ber_qZ=\det_qA\cdot\det_qV. \]
\end{pro}
\proof First, we give the formula for the quantum determinant $\det_qA$ and $\det_qV$. Let $\Phie=R_{e(1, 1, \ldots ,1)}$ and $\Phio=\Ro_{(s)}$ in the notation of the proof of Lemma \ref{lem41}. Then one can show that $\Phi_e$ is the projection on $\Lambda_0^{\vee}$ and $\Phi_o$ is the projection on $S_1^{\vee s}$ (cf. \cite{gur}). Let $S_e$ be the matrix such that $\Re_{im}^{jn}\Se^{mk}_{nl}=\Se_{im}^{jn}\Re^{mk}_{nl}=\delta_{i}^{k}\delta^j_l$. Let $\Fe_i^j:=\Se_{li}^{lj}$. We define the matrix $\So$ analogously and set $\Go_i^j:=\So_{il}^{jl}$.

\begin{lem}\label{lem44} Let ${\cal I}(r,d)=\{I=(i_1, i_2, \cdots , i_r)| i_j\in \{1, 2, \cdots d\}\}$ be the set of indices. For $I=(i_1, \cdots i_r), J=(j_1, \cdots j_r)$ let $\Fe_I^J=\Fe_{i_1}^{j_1}\Fe_{i_2}^{j_2}\cdots \Fe_{i_r}^{j_r}$. Then the quantum determinant for $A$ can be given by 
\[ \det_qA=q^{\frac{r(r+1)}{2}}\Fe_J^I {\Phi_e}_K^JA_I^K,\]
where the indices run in $\cal I$. \end{lem}
\proof
Let us consider the bialgebra $E_0$ generated by entries of $A$ and its Hopf envelope $H_0$. 
We know that $\Lambda_0^r$ defines the quantum determinant for $A$. The coaction $\delta:\Lambda_0^r\longrightarrow\Lambda_0^r\otimes H_0$ induces the coalgebra homomorphism $\overline{\delta}: {(\Lambda_0^r)}^*\otimes \Lambda_0^r\longrightarrow H_0$, whose image is the subspace of $H_0$, generated by the quantum determinant. ${(\Lambda_0^r)}^*\otimes\Lambda_0^r$ is one-dimensional, hence it is sufficient to find a non-zero element of this space and consider its image in $H_0$. 

Since $R_e$ is closed, we can define the map $\varphi_1:{\Bbb K}\stackrel{\db}{\longrightarrow}V\otimes V^*\stackrel{\Se}{\longrightarrow} V^*\otimes V,$ $\db$ is defined as on page \pageref{differential}. We have $\varphi_1(1_{{\Bbb K}})=\Fe_j^i\xi^j\otimes x_i$. More generally, we define the map $\varphi_r:{\Bbb K}\longrightarrow V_0^{\otimes r}\otimes{(V_0^{\otimes r})}^*\longrightarrow {(V_0^{\otimes r})}^*\otimes V_0^{\otimes r}$, $\varphi_r(1_{\kk})=\Fe^I_J\xi^J\ot x_I$.
Therefore we get the map $(\Phie^*\otimes\Phie)\circ\varphi_r:{\Bbb K}\rightarrow{(\Lambda_0^r)}^*\otimes\Lambda_0^r,$ whence, $\overline{\delta}_o(\Phie^*\otimes \Phie)\circ\varphi_r:{\Bbb K}\longrightarrow H$. The image of $1_{{\Bbb K}}$ under $\overline{\delta}$ is then $\Fe_J^I{\Phi_e}_K^JA_L^K\Phie_I^L$. Since $\Phie$ and $A$ commute and since $\Phie$ is a projection, this element is equal to $\Fe_J^I{\Phie}_K^JA_I^K$. To show that this element is the quantum determinant we are looking for, it is enough to show that it is not zero. We have $\varepsilon (\Fe_J^I{\Phie}_K^JA_I^K)=\Fe_J^I{\Phie}_I^J=q^{-\frac{r(r+1)}{2}}$ where $\varepsilon$ is the counit of $H_0$. The last equation is due to the fact that $\Fe_j^i{\Re_i^j}_l^k=\delta_l^k$ and $\Fe_i^i=\frac{1-q^{-r}}{q-1}$ (cf. \cite{ph97b}).\eee

In the same way we get the formula for the quantum determinant for the quantum matrix $V$:
\[\det_qV={q}^{\frac{s(s+1)}{2}}\Go_J^I\overline{\Phio}_K^JV_I^K,\]
where $\overline{\Phio}_{k_1k_2\cdots k_n}^{j1j2\ldots j_n}:=\Phio_{k_n\ldots k_2k_1}^{j_n\ldots j_2j_1}$. The order of indices of $\overline{\Phio}$ is reversed because of the definition of the pairing on page \pageref{pairing}. However, using the relation in \rref{ttrelation} and the fact that $\Phio$ commutes with all $\Ro_i, i=1,2,\ldots,n$, we deduce that $\det_qV={q}^{\frac{s(s+1)}{2}}\Go_J^I\Phio_K^JV_I^K$.

\begin{lem}\label{lem45} $\det_qA \mbox{ and } \det_qV$ commute.\end{lem}
\proof
Let $T$ be the matrix such that $T_{iq}^{jp}P_{pl}^{qk}=\delta_l^k\delta_l^j$ then, according to \rref{ztrelation}, we have
\[A_j^iV_k^l=T_{jp}^{lq}V_n^pA_q^mP_{mk}^{ni}.\]
Hence, there exist matrices $T_{KL}^{IJ}$ and $P_{KL}^{IJ}$ such that
\begin{eqnarray*}
  T_{IQ}^{JP}P_{PL}^{QK} = \delta_I^K\delta_L^J 
& \mbox{and}& 
A^K_IV^P_M=T^{PR}_{IS}V^S_UA^V_RP^{UK}_{VM}.  \end{eqnarray*}
 The equations in \rref{eqforp} imply
\[ \begin{array}{c@{\quad}c}
\Fe_J^IT_{IS}^{PR}=\Fe_I^R {P^{-1}}_{SJ}^{IP}, & \Phie^I_JP^{MJ}_{KN} = P^{MI}_{JN}\Phi^J_K,\\ 
{P^{-1}}_{SJ}^{IP}\Go_P^N=\Go_S^PT_{JP}^{NI}, & \Phio^I_JP^{JM}_{NK} =P^{IM}_{NJ}\Phio^J_K. \end{array} \]
Hence \[ \Fe_J^IT_{IS}^{PR}\Go_P^N=\Fe_I^R\Go_S^PT_{JP}^{NI},\quad
\Phie_K^JP^{UK}_{VM}\Phio_N^M=\Phie^K_VP_{KN}^{MJ}\Phio_M^U.\]
Therefore, we have:
\begin{eqnarray*}
\Fe_J^I\Phie_K^JA_I^K\cdot\Go_N^M\Phio_P^NV_M^P & = & \Fe_J^I\Phie_K^JT^{PR}_{IS}V_U^SA_R^VP^{UK}_{VM}\Go^N_P\Phio_N^M\\
& = & \Phie_K^J\Fe_I^R\Go_S^PT_{JP}^{NI}P^{UK}_{VM}\Phio_N^M\ V_U^SA_R^V\\
& = & \Phie^K_V\Fe_I^R\Go_S^PT_{JP}^{NI}P_{KN}^{MJ}\Phio_M^U\ V_U^SA_R^V\\
& = &\Phio_M^U\Go_S^MV_U^S\cdot \Phie^I_V\Fe_I^RA_R^V
\end{eqnarray*}
The lemma is proved.\eee

We now prove the formula for the Berezinian. Let $\alpha_J^I$ and $\beta_L^K, I, J \in {\cal I}(r,m)$ and $K, L \in {\cal I}(s,n)$ be two sets of indices, such that $\alpha^IX_I\beta_J\Xi^J$ represents the unique cohomology class of $H(K_0^{\bullet\bullet})$. Thus 
\bba\label{hkber} &\delta ( \alpha^IX_I\beta_J\Xi^J)=\alpha^IX_I\beta_J\Xi^J\otimes \Ber  \mbox{  (modulo  Im} \d\otimes H_0).& \eea

On the other hand, $\alpha^IX_I$ (resp. $\beta_J\Xi^J$) represents $\Lambda_0^r$ (resp. $S_1^{\vee s}$) and 
\[\delta ( \alpha^IX_I)=\alpha^IX_I\otimes \det_qA+X_K\otimes Z^K\]
where $K$ runs in the set ${\cal I}(r, d)\backslash {\cal I}(r ,m)$ and $Z^K$ are certain elements of $H$. A similar equation holds for $\det_qV$:
\[\delta ( \beta_J\Xi^J)=\beta_J\Xi^J\otimes\det_qV+\Xi^L\otimes T_L , L\in {\cal I}(s, d)\backslash {\cal I}(s, n)\]
Thus, \[\delta ( \alpha^IX_I\beta_J\Xi^J)=\alpha^IX_I\beta_J\Xi^J\otimes\det_qA\det_qV+X_K\xi^L\otimes Z^KT_L\]
where $(K, L)$ runs in ${\cal I}(r, d)\times {\cal I}(s, d)\setminus {\cal I}(r, m) \times {\cal I}(s, n)$. We  notice that the monomials $X_I\Xi^J, (i, j)\in 	{\cal I}(r, m) \times {\cal I}(s, n)$ and $X_K \Xi^L, (K, L) \in {\cal I}(r, d) \times {\cal I}(s, d)\backslash {\cal I}(r, m)$ $ \times {\cal I}(s, n)$ form two subspaces of $\Lambda^r\otimes S^{\vee s}$, that intersect each other at 0. According to \rref{hkber}, one should have 
\[X_K\xi^L\otimes Z^Kt_L\equiv 0 \quad(\mbox{ modulo } \Im \d\otimes H).\]
Therefore, $\Ber_qZ=\det_qA\det_qV$, the proof of Proposition \ref{pro43} is completed.\eee
\begin{cor}
If $Z$ is invertible then so are $A$ and $D$. \end{cor}
Like the case of multiparameter deformations, we can easily derive from the above that $D-CA^{-1}B$ and $A-BD^{-1}C$ are also invertible quantum matrices and an analog of Lemma \ref{lemisdet} follows from Lemma \ref{lem44}. Hence
\begin{eqnarray*}
\Ber_qZ & =& \det_qA\det_q{(D-CA^{-1}B)}^{-1}\\
&=&\det_qD^{-1}\det_q(A-CD^{-1}B).\end{eqnarray*}
Thus, Theorem \ref{firstimp} holds for matrix quantum super groups of type $(r, s)$. Its converse, i.e. Theorem \ref{edl6}, also holds. In fact, everything can be done along the same line as in \S\ref{sec6}. There is no difference for the proof of Lemma \ref{lem1}, from which Corollary \ref{cor2} follows. The same situation is for the proofs of Lemma \ref{lem4} and Corollary \ref{cor5} from which Theorem \ref{edl6} follows. Thus we have
\begin{edl}\label{lastdl}
Let $Z$ be a quantum super matrix in a standard format, associated to the Hecke sum of an even Hecke operator and an odd Hecke operator. Then $Z$ is invertible iff its even block matrices $A$ and $D$ are invertible. The quantum Berezinian for $Z$ can be defined as in the classical case and is multiplicative.\end{edl}

\end{document}